\documentclass[11pt]{article}

\usepackage{IEEEtrantools}
\usepackage{latexsym}
\usepackage{amsfonts}
\usepackage{amssymb}
\usepackage{amsmath,amsfonts,amsthm}
\usepackage{mathrsfs}
\usepackage{enumerate}
\usepackage{dsfont}
\usepackage{stix}
\usepackage{bm}
\usepackage{multirow}
\usepackage{tikz}
\usepackage{tikz-cd}
\usetikzlibrary{arrows,decorations.markings}
\usepackage{bookmark}
\usepackage{hyperref}
\hypersetup{
     colorlinks   = true,
     citecolor    = blue
}

\newcommand{\bbbt}{\mathbb{T}}
\newcommand{\scrt}{\mathscr{T}}

\newcommand{\be}{\begin{equation}}
\newcommand{\ee}{\end{equation}}
\newcommand{\bea}{\begin{eqnarray}}
\newcommand{\eea}{\end{eqnarray}}
\newcommand{\bean}{\begin{eqnarray*}}
\newcommand{\eean}{\end{eqnarray*}}
\newcommand{\brray}{\begin{array}}
\newcommand{\erray}{\end{array}}
\newcommand{\biearray}{\begin{IEEEarray}{rCl}}
\newcommand{\eiearray}{\end{IEEEarray}}

\newtheorem{dfn}{Definition}[section]
\newtheorem{thm}[dfn]{Theorem}
\newtheorem{lmma}[dfn]{Lemma}
\newtheorem{ppsn}[dfn]{Proposition}
\newtheorem{crlre}[dfn]{Corollary}
\newtheorem{xmpl}[dfn]{Example}
\newtheorem{rmrk}[dfn]{Remark}

\newcommand{\bdfn}{\begin{dfn}\rm}
\newcommand{\bthm}{\begin{thm}}
\newcommand{\blmma}{\begin{lmma}}
\newcommand{\bppsn}{\begin{ppsn}}
\newcommand{\bcrlre}{\begin{crlre}}
\newcommand{\bxmpl}{\begin{xmpl}}
\newcommand{\brmrk}{\begin{rmrk}\rm}

\newcommand{\edfn}{\end{dfn}}
\newcommand{\ethm}{\end{thm}}
\newcommand{\elmma}{\end{lmma}}
\newcommand{\eppsn}{\end{ppsn}}
\newcommand{\ecrlre}{\end{crlre}}
\newcommand{\exmpl}{\end{xmpl}}
\newcommand{\ermrk}{\end{rmrk}}

\newcommand{\bbc}{\mathbb{C}}
\newcommand{\bbz}{\mathbb{Z}}
\newcommand{\bbn}{\mathbb{N}}

\newcommand{\clh}{\mathcal{H}}

\newcommand{\clk}{\mathcal{K}}

\newcommand{\cls}{\mathcal{S}}

\def \bbt {\mbox{\boldmath $t$}}

\newcommand{\prf}{\noindent{\it Proof\/}: }

\newcommand{\one}{{1\!\!1}}

\def \qed { \mbox{}\hfill
$\Box$\vspace{1ex}}

\begin{document}

%%%%%%%%%%%%%%%%%%%%%%%%%%%%%%%%%
%%%%%%%%%%%%%%%%%%%%%%%%%%%%%%%%%

\author{\sc{Bipul Saurabh}}
\title{ $C(SO_q(2n+1)/SO_q(2n-1))$ as iterated torsioned quantum double suspensions of   $C(\bbbt)$ }
\maketitle

%%%%%%%%%%%%%%%%%%%%%%%%%%%%%%%%%%
%%%%%%%%%%%   ABSTRACT    %%%%%%%%%%%%%%%%  
%%%%%%%%%%%%%%%%%%%%%%%%%%%%%%%%%%

\begin{abstract}
Let $A$ be a unital $C^*$-algebra, and let $\Sigma^2_m A$ denote the $m$-torsioned quantum double suspension of $A$. 
For $q \in (0,1)$ and $n \geq 1$, we prove that the $C^*$-algebra corresponding to the quotient space
$SO_q(2n+1)/SO_q(2n-1)$ is isomorphic to 
$
\Sigma^{2(n-1)} \, \Sigma^2_2 \, \Sigma^{2(n-1)} C(\mathbb{T})$.
It follows as a consequence that these spaces are independent of the deformation parameter
$q$. 

\end{abstract}
\bigskip

{\bf AMS Subject Classification No.:} {\large 19}K{\large 33},  {\large 46}L{\large 80},
 {\large 58}B{\large 34}. \\

{\bf Keywords.}  $C^*$-extension,  homogeneous extension,  corona factorization property, $m$-torsioned quantum double suspension.
\bigskip

%%%%%%%%%%%%%%%%%%%%%%%%%%%%%%%%%%
%%%%%%%%%   INTRODUCTION    %%%%%%%%%%%%%%%
%%%%%%%%%%%%%%%%%%%%%%%%%%%%%%%%%%

\section{Introduction}

Hong and Szymanski introduced the notion of quantum double suspension (QDS) for unital 
$C^*$-algebras in \cite{HonSzy-2002aa}. This construction allows one to pass from classical compact spaces to noncommutative ones in a systematic way.  When a noncommutative space is obtained by iteratively applying the QDS to a compact space $X$, many of its topological and geometric properties can often be understood in terms of those of $X$. For instance, it simplifies the computation of $K$-groups and produces examples of noncommutative geometries from the geometries of the classical space (see \cite{ChaSun-2011ab}). 
Similar to QDS, for any unital $C^*$-algebra $A$,  the $m$-torsioned quantum double suspension, denoted by $\Sigma^2_m A$, is defined in   \cite{BhuBisSau-2024aa}. For $m=1$, this construction reduces to the usual QDS. As in the case of QDS, it is desirable to realize a given $C^*$-algebra as an iterated $m$-torsioned quantum double suspension of $C(X)$ for some compact space $X$. In \cite{BhuBisSau-2024aa}, it was shown that the $C^*$-algebra $C(SO_q(3))$ is isomorphic to $\Sigma^2_2 C(\mathbb{T})$.  
The purpose of the present article is to generalize this result and show that the $C^*$-algebra associated with the quantum homogeneous space $SO_q(2n+1)/SO_q(2n-1)$ can be obtained from $C(\mathbb{T})$ by successively applying the QDS and the $2$-torsioned QDS in a suitable order. More precisely, we establish an isomorphism
between $
C\!\left(SO_q(2n+1)/SO_q(2n-1)\right)$ 
 and 
$\Sigma^{2(n-1)} \, \Sigma^2_2 \, \Sigma^{2(n-1)} C(\mathbb{T})$.
This description immediately implies  that the topological type of the underlying $C^*$-algebra is independent of the deformation parameter $q$. 

We now briefly outline the main idea. For a nuclear $C^*$-algebra $A$ and a finite-dimensional compact metric space $Y$ (that is, a closed subset of $\cls^n$ for some $n \in \bbn$), Pimsner, Popa, and Voiculescu~\cite{PimPopVoi-1979aa} introduced a group $\operatorname{Ext}_{\mathrm{PPV}}(Y, A)$ consisting of strongly unitary equivalence classes of homogeneous extensions of $A$ by $C(Y) \otimes \clk$. 
An important feature of this group, which distinguishes it from the usual $Ext$ group, is that any two elements representing the same class must have isomorphic middle $C^*$-algebras. 
This property plays a crucial role in the present work. 
We first establish an isomorphism between $\operatorname{Ext}_{\mathrm{PPV}}(Y, A)$ and $\operatorname{Ext}_{\mathrm{PPV}}(Y, \Sigma^2_m A)$, under  mild assumptions on $A$ and $Y$. 
Taking $Y = \bbbt$, we then show that $K_0$ serves as a complete invariant for the class of $C^*$-algebras arising as middle algebras of homogeneous extensions of $\Sigma^2_2 \Sigma^{2(n-1)} C(\bbbt)$ by $C(\bbbt) \otimes \clk$. 
The result follows by comparing the $K$-groups, and then extending it  using the main result of   \cite{Sau-2019aa}.

Our article is organized as follows.  In Section~2, we recall the construction of $m$-torsioned quantum double suspensions, and the quotient space $SO_q(2n+1)/SO_q(2n-1)$, and and some related results.  In Section~3, we establish the isomorphism between  $\operatorname{Ext}_{\mathrm{PPV}}(\bbbt, A)$ and $\operatorname{Ext}_{\mathrm{PPV}}(\bbbt, \Sigma^2_mA)$. Using this, we prove that $C\!\left(SO_q(2n+1)/SO_q(2n-1)\right)$ is isomorphic to $\Sigma^{2(n-1)}\Sigma^2_2 \Sigma^{2(n-1)} C(\mathbb{T})$.

We now fix some notation. Throughout, the symbol $q$ denotes a real number in the interval $(0,1)$.
The standard orthonormal bases of the Hilbert spaces $\ell^{2}(\bbn)$ and $\ell^{2}(\bbz)$ are denoted by
$\{e_n : n \in \bbn\}$ and $\{e_n : n \in \bbz\}$, respectively.
The left shift operator on both $\ell^{2}(\bbn)$ and $\ell^{2}(\bbz)$ will be denoted by the same symbol $S$.
For $m<0$, the expression $(S^{*})^{m}$ denotes the operator $S^{-m}$.
Let $p_{ji}$ denote the rank-one operator sending $e_i$ to $e_j$.
We write $p_{ii}$ simply as $p_i$, and $p_{00}$ as $p$.
The symbol $p_{<m}$ denotes the projection $\sum_{i=0}^{m-1} p_i$.
We write $\mathcal{L}(\clh)$ and $\clk(\clh)$ for the sets of all bounded linear operators and all compact operators
on a Hilbert space $\clh$, respectively, and denote by $\clk$ the $C^{*}$-algebra of compact operators.
For a $C^{*}$-algebra $A$, we use $M(A)$ and $Q(A)$ to denote its multiplier algebra and corona algebra, respectively.
The map $\pi$ denotes the canonical homomorphism from $M(A)$ onto $Q(A)$, and for $a \in M(A)$,
the symbol $[a]$ denotes its image under $\pi$. For a group $G$,  let $G^{\mathrm{Tor}}$ denote the torsion subgroup of $G$.

    \section{Preliminaries}
 We begin by recalling the definition and certain known results concerning the $m$-torsioned quantum double suspension and the quotient space $SO_q(2n+1)/SO_q(2n-1)$, which will be used in the later sections.

    \subsection{$m$-torsioned Quantum double suspension}
    \bdfn \label{QDS}  Let $A$ be a unital $C^*$-algebra. For $m \in \bbz\setminus \{0\}$, we define its $m$-torsioned quantum double suspension  as the unital $C^*$-algebra $\Sigma_{m}^{2}A$ for which there exists an essential extension 
    \[
    0 \rightarrow A \otimes \clk \rightarrow \Sigma_{m}^{2}A \rightarrow C(\bbbt) \rightarrow 0 
    \]
    such that the corresponding Busby invariant $\beta: C(\bbbt) \rightarrow Q(A \otimes \clk)$ mapping $\beta(t)=[1 \otimes (S^*)^{m}]$. Equivalently, $\Sigma_{m}^{2}A$ can be defined as the $C^*$-subalgebra of $\scrt \otimes A$ generated by $A \otimes \clk$ and $1 \otimes (S^*)^{m}$. By $\Sigma_m^{2n}A$ we mean the $n$-fold iteration of $\Sigma_m^{2}$ applied to $A$, that is, $\underbrace{\Sigma_m^{2}\Sigma_m^{2}\cdots \Sigma_m^{2}A}_{n-\mbox{copies}}$.
    \edfn

  The following proposition establishes the  universal property of $m$-torsioned quantum double suspension.
    \bppsn  \rm{(\cite{BhuBisSau-2024aa})} \label{universal} Let $\phi:A \longrightarrow B$ be a  homomorphism with $\phi(1)=P$. Let $T \in B$ be an isometry with the defect projection $P$ and  $\nu:M_m(\bbc) \rightarrow B$ be a $*$-homomorphism satisfying 
    \begin{IEEEeqnarray}{rCl} \label{eq1}
    	\nu(1)=P \mbox{ and } \quad \nu(D)\phi(a)=\phi(a)\nu(D), \, \mbox { for all } a \in A \mbox { and  } D \in M_m(\bbc). 
    \end{IEEEeqnarray} 
    Then there exists a unique $*$-homomorphism $\Sigma_m^2 (\phi,\nu,T):\Sigma_m^2A \rightarrow B$ such that 
    $$ \Sigma_m^2 (\phi, \nu,T)(a \otimes p_{ij})=\phi(a) \nu(p_{ij}),  \mbox{ for } 0\leq i,j \leq m-1 $$
    and 
    $$ \Sigma_m^2 (\phi,\nu,T)(1 \otimes (S^*)^m)=T.$$
    Conversely, let  $\psi:\Sigma_m^2A \rightarrow B$ be any unital $*$-homomorphism. Define $T=\psi(1\otimes (S^*)^m)$ and $P=1-TT^*$. Then there exist homomorphisms 
    $ \phi:A \longrightarrow B$ and $\nu:M_m(\bbc) \rightarrow B$, and an isometry $T$  satisfying equation (\ref{eq1}) such that 
    \[
    \psi =\Sigma_m^2 (\phi, \nu,T).
    \]
    \eppsn
    
  \bppsn \label{ideals}
  Any closed ideal $J$ of $\Sigma_m^2 A$ is either of the form $I \otimes \clk$, where $I$ is a closed ideal in  $A$, or it contains $A \otimes \clk$. 
  \eppsn 
  \prf Given a representation $\rho$ of $A$ acting on the Hilbert space $\clh$, one can define a representation $\Sigma^2_m \rho$ of $\Sigma^2_m A$ actiong on $\clh \otimes \ell^2(\bbn)$ as follows (see \cite{BhuBisSau-2024aa}).
  $$ \Sigma^2 \rho\,(1\otimes S^*)=1\otimes S^* \, \, \mbox { and   } \, \, \Sigma_m^2 \rho\,(a\otimes p)=\rho(a)\otimes p, \mbox{ for all } a \in A. $$ 
It is not difficult to verify  that $\ker \Sigma^2_m \rho=\ker \rho \otimes \clk$.  Using this and Theorem $5.7$ in \cite{BhuBisSau-2024aa}, the claim follows.
  \qed 
    
    \bppsn \label{ideals 1}
    Let $\phi: \Sigma_m^2A \rightarrow B$ be a $*$-homomorphism such that $\phi(a \otimes k)\neq 0$ for any $a \in A$ and $k \in \clk$. Then $\phi$ is injective. 
    \eppsn 
    \prf  Since $\phi(a \otimes k)\neq 0$ for any $a \in A$ and $k \in \clk$, it follows from Proposition \ref{ideals} that the only possibility is $\ker \phi=0$. 
    \qed 
    \bthm  \rm{(\cite{BhuBisSau-2024aa})} \label{K-groups for m-QDS}
    Let $K_0(A)$ and $K_1(A)$ be finitely generated abelian groups with $\bbz$-linearly independent 
    generators $\left\{\left[P_i\right]\right\}_{i = 1}^{r}$ and  $\left\{\left[U_i\right]\right\}_{i = 1}^s$,  respectively.
  Assume that $[1] =[P_1]$.  Then
    $K_0(\Sigma_m^2A)$ is isomorphic to $K_0(A)\oplus \bbz/m\bbz$.  
    The  generators $[1]$,$\left\{[P_i\otimes p]\right\}_{i=2}^{r}$  generate the subgroup $K_0(A)$ of $K_0(\Sigma_m^2A)$ and $[1\otimes p]$ generates
    the component $\bbz/m\bbz$ of   $K_0(\Sigma_m^2A)$. Moreover, the group $K_1(\Sigma_m^2A)$ 
    is isomorphic to $K_1(A)$ with generators $\left\{\left[U_i\otimes p+1-1\otimes p\right]\right\}_{i = 1}^s$.
    \ethm
    
    \subsection{The quotient space $SO_q(2n+1)/SO_q(2n-1)$}
    Let \[ 
    B_k^{2n+1}(q)= 
    \begin{cases}
    	C(\bbbt) & \mbox{ if } k=1, \cr
    	\phi_{\omega_k}(C(SO_q(2n+1)/SO_q(2n-1))) & \mbox{ if } 1<k \leq 2n. \cr
    \end{cases}
    \]
We now recall some results from \cite{BhuBisSau-2024aa}, which are needed to prove our main claim.
   \blmma   \rm{(\cite{BhuBisSau-2024aa})}\label{lemma-homogeneous}
   For $1 < k \leq 2n$, the short exact sequence 
   \[
   \chi_{k}: \quad 0\longrightarrow C(\bbbt) \otimes \clk \xrightarrow{i} B_k^{2n+1}(q) \xrightarrow{\rho_{k+1}}  B_{k-1}^{2n+1}(q) \longrightarrow  0.
   \]
   is a unital homogeneous extension of $B_{k-1}^{2n+1}$ by 
   $C(\bbbt)\otimes \clk$. 
   \elmma             
   
    \bthm  \rm{(\cite{BhuBisSau-2024aa})} \label{K-groups-B}
    For $1 \leq k \leq 2n$, define $u_k= t \otimes p^{\otimes (k-1)} +1-1\otimes  p^{\otimes (k-1)}$. 
    Then one has 
    \begin{IEEEeqnarray*}{rCl}
    	K_0(B_{k}^{2n+1})&=&\begin{cases}
    		\langle[1]\rangle \cong \bbz, & \mbox{ if } 1 \leq  k \leq n, \cr 
    		\langle[1]\rangle \oplus \langle[1 \otimes p^{\otimes (n) }\otimes 1^{\otimes (k-n-1)}]\rangle  \cong \bbz \oplus \bbz/2\bbz, & \mbox{ if } n+1 \leq k \leq 2n, \cr 
    	\end{cases} \\
    	K_1(B_{k}^{2n+1})&=&\langle [u_k] \rangle\cong \bbz.
    \end{IEEEeqnarray*}
    \ethm
    \subsection{$C^*$-algebra extensions}
 Let $A$,  and $B$ be unital separable nuclear $C^*$-algebras. 
  Two elements $a,b \in Q(B \otimes \mathcal{K})$ are said to be \emph{strongly unitarily equivalent} 
  if there exists a unitary $U \in M(B \otimes \mathcal{K})$ such that
  \[
  [U]\, a \, [U^*] = b,
  \]
  and we write $a \sim_{su} b$.
  Two $C^*$-extensions $\alpha, \beta : A \to Q(B \otimes \mathcal{K})$ are said to be 
  \emph{strongly unitarily equivalent}, denoted by $\alpha\sim_{su}\beta $,  if there exists a unitary 
  $U \in M(B \otimes \mathcal{K})$ such that
  \[, 
  [U]\, \alpha(a) \, [U^*] = \beta(a)
  \quad \text{for all } a \in A.
  \]
  An element $a$ in a $C^*$-algebra $A$ is called \emph{norm-full} 
  if it is not contained in any proper closed ideal of $B$.
  An extension
  \[
  \tau : A \to Q(B \otimes \mathcal{K})
  \]
  is said to be \emph{norm-full} if for every nonzero element $a \in A$, 
  the element $\tau(a)$ is norm-full in $Q(B \otimes \mathcal{K})$.
   \bdfn
  	Let $B$ be a separable stable $C^*$-algebra. 
  	Then $B$ is said to have the \emph{corona factorization property} 
  	if every norm-full projection in $M(B)$ is Murray-von Neumann equivalent 
  	to the unit element of $M(B)$.
  \edfn
   It is easy to see that if a $C^*$-algebra $B$ has the corona factorization property, 
  then any norm-full projection in $Q(B)$ is Murray-von Neumann equivalent 
  to the unit of $Q(B)$. 
  Furthermore, one can show that for a finite-dimensional compact metric space $Y$, 
  the algebra $C(Y) \otimes \mathcal{K}$ has the corona factorization property 
  (see \cite{PimPopVoi-1979aa}).

Suppose that $Y$ is a finite-dimensional compact metric space. 
In other words, $Y$ is homeomorphic to a closed subset of the Euclidean sphere $\cls^n$ for some $n \in \mathbb{N}$.
Let
  \[
  M(Y) := M(C(Y) \otimes \mathcal{K}), 
  \quad 
  Q(Y) := M(C(Y) \otimes \mathcal{K}) / (C(Y) \otimes \mathcal{K}),
  \]
  and let
$Q := \mathcal{L}(\mathcal{H}) / \mathcal{K}(\mathcal{H})$
  denote the Calkin algebra.
   It is easy to see that $M(Y)$ can be identified with the set of all 
  $*$-strongly continuous functions from $Y$ to $\mathcal{L}(\mathcal{H})$.
 An extension $\tau$ of $A$ by $C(Y) \otimes \mathcal{K}$ is said to be 
  \emph{homogeneous} if for every $y \in Y$, the map
  \[
  \mathrm{ev}_y \circ \tau : A \rightarrow Q
  \]
  is injective, where $\mathrm{ev}_y : Q(Y) \rightarrow Q$ denotes the evaluation map at $y$.
  Let $\mathrm{Ext}_{\mathrm{PPV}}(Y,A)$ denote the set of unitary equivalence classes 
  of unital homogeneous extensions of $A$ by $C(Y) \otimes \mathcal{K}$. 
  For a nuclear $C^*$-algebra $A$, 
  Pimsner, Popa, and Voiculescu \cite{PimPopVoi-1979aa} showed that 
$\mathrm{Ext}_{\mathrm{PPV}}(Y,A)$  is a group. 
  We denote the equivalence class of an extension $\tau$ in 
$\mathrm{Ext}_{\mathrm{PPV}}(Y,A)$  by $[\tau]_{su}$. Denote by $[\tau]$ the stable equivalence class of an extension $\tau$ in the group $\mathrm{Ext}(A,C(Y))=KK^1(A,C(Y))$. 
    \bppsn \rm{(\cite{Sau-2019aa})}  \label{injectivekk}
  Let $A$ be a unital separable nuclear $C^*$-algebra satisfying the Universal Coefficient Theorem. 
  Suppose that $Y$ is a finite-dimensional compact metric space. 
  Then the map
  \[
  i :\mathrm{Ext}_{\mathrm{PPV}}(Y,A)  \longrightarrow KK^1(A,C(Y)), 
  \qquad
  [\tau]_{su} \longmapsto [\tau]
  \]
  is an injective homomorphism.
  \eppsn
  
    \section{Isomorphism between $\operatorname{Ext}_{\mathrm{PPV}}(Y, A)$ and $\operatorname{Ext}_{\mathrm{PPV}}(Y, \Sigma^2_mA)$}
  In this section, we begin by proving that, under suitable assumptions on the space $Y$ and on a unital, nuclear, and separable 
 $C^*$-algebra $A$, the groups $\operatorname{Ext}_{\mathrm{PPV}}(Y, A)$ and $\operatorname{Ext}_{\mathrm{PPV}}(Y, \Sigma^2_mA)$ are isomorphic. This is then used to establish the main result of the paper. 
   Without loss of generality, we shall assume that the Hilbert space $\clh$ is $\ell^2(\bbn)$.  Throughout this section,  $Y$ is a finite-dimensional compact metric space.

\bdfn 
Let $m \geq 0$ and $\ell \geq 2$.    
A $\ast$-homomorphism   
\[
\nu: M_m(\bbc) \longrightarrow Q\!\left(C(Y)\otimes \clk^{\otimes \ell}\right), 
\quad  \nu(1)= [1^{\otimes \ell} \otimes p_{<m}]
\]
is called an $m$-torsion system of length $\ell$. Define 
\[
\Delta_m^{\ell}=\left\{\nu: M_m(\bbc) \to Q\!\left(C(Y)\otimes \clk^{\otimes \ell}\right) 
\ \middle|\ \nu(1)= [1^{\otimes \ell} \otimes p_{<m}] \right\}.
\]
\edfn 

\brmrk Note that any $\nu \in \Delta_m^{\ell}$ is an essential norm full  extension of $M_m(\bbc)$ by $Q(C(Y)\otimes \clk$. Its stable equivalence class  $[\nu]$ is an element  of the group  $\operatorname{Ext}(M_m(\mathbb{C}), C(Y))$. 
\ermrk

Let $\tau$ be a unital homogeneous extension of $A$ by $C(Y)\otimes \clk(\clh)$.  
We define
\[
\tilde{\tau} : A \longrightarrow Q\big(C(Y)\otimes \clk(\clh) \otimes \clk(\clh)\big), 
\qquad 
\tilde{\tau}(a) =  [\tau(a)_y\otimes p_{<m}]_{y\in Y}.
\]
Take $\nu \in \Delta_m^2$. Define the extension $\Sigma_m^2 (\tau, \nu):=\Sigma^2_m(\tilde{\tau}, \nu, [1\otimes 1  \otimes S^m])$. Equivalently, one may write $\Sigma_m^2 (\tau, \nu)$ as the $*$-homomorphism given by:
\begin{equation}\label{sigmatau}
	\Sigma_m^2 (\tau, \nu) : \Sigma_m^2 A \;\longrightarrow\; Q\big(C(Y)\otimes \clk(\clh)\otimes \clk(\clh)\big)
\end{equation}
such that
\[
	\Sigma_m^2 (\tau, \nu)(a \otimes p_{<m}) = \tilde{\tau}(a) = [\tau(a)_y \otimes \mathcal{P}_{<m}]_{y \in Y},
\quad 
	\Sigma_m^2 (\tau, \nu)(1 \otimes S^m) = [1 \otimes S^m]_{y \in Y} \]
and 
\[ 
	\Sigma_m^2 (\tau, \nu)(1 \otimes p_{ij})= [1 \otimes \nu(p_{ij})_y]_{y \in Y}, \mbox{ for } 0\leq i,j \leq m-1.
\]
Since $\tau$ is homogeneous, one can check that  
$$\mbox{ev}_y \circ \Sigma_m^2 (\tau, \nu)(a \otimes k)\neq 0$$ 
 for any $a \in A$ and $k \in \clk$. As a consequence of Propostion \ref{ideals 1},  it follows from  that the extension 
 $\Sigma_m^2 (\tau, \nu)$ is  homogeneous. If $\nu=0$, then we denote  $\Sigma_m^2 (\tau, \nu)$ by  $\Sigma_m^2\tau$. If $m=1$, then we write $\Sigma_m^2 (\tau, \nu)$ and $\Sigma_m^2\tau$ by 
 $\Sigma^2 (\tau, \nu)$ and $\Sigma^2\tau$. Observe that  
the extension  $\Sigma_m^2 \tau$ is nothing but the restriction of  the map  	$\Sigma^2 \tau$ to the subalgebra  $\Sigma_m^2A$. 
It follows from Lemma $2.8$ in \cite{Sau-2019aa} that 
the map 
$$\beta: \operatorname{Ext}_{\mathrm{PPV}}(Y,  A)  \longrightarrow \operatorname{Ext}_{\mathrm{PPV}}(Y,  \Sigma^2A); \quad [\tau]_{su} \;\longmapsto\; [\Sigma^2 \tau]_{su}$$ 
is an isomorphism.
Let 
$$ i : \operatorname{Ext}_{\mathrm{PPV}}(Y,  \Sigma^2A) \rightarrow \operatorname{Ext}_{\mathrm{PPV}}(Y,  \Sigma_m^2A); \quad [\Sigma^2 \tau]_{su} \mapsto [\Sigma_m^2 \tau]_{su}.$$
It is easy to see that $i$ is well-defined homomorphism.  Therefore, by composing it with $\beta$, we get a homomorphism 
$$ \Psi_m: \operatorname{Ext}_{\mathrm{PPV}}(Y,  A)  \longrightarrow \operatorname{Ext}_{\mathrm{PPV}}(Y,  \Sigma_m^2A); \qquad \phi=i \circ \beta.$$ 
\blmma \label{injective hom}
Let $m \geq 1$, and  let $A$ be a unital, nuclear, separable 
$C^*$-algebra. Let $Y$ be a finite-dimensional compact metric space.  Then the map
\[
\Psi_m: \operatorname{Ext}_{\mathrm{PPV}}(Y, A) \to \operatorname{Ext}_{\mathrm{PPV}}(Y, \Sigma_m^2 A), \quad [\tau]_u \mapsto [\Sigma_m^2 \tau]_u
\]
is an  injective homomorphism. 
\elmma
\prf  
Assume that $\Psi_m([\tau]_{su})=[\Sigma_m^2\tau]_{su}=0$.  
Then there exists a $*$-homomorphism $\overline{\Sigma^2_m\tau}:  \Sigma_m^2A\rightarrow   M(C(Y) \otimes \clk(\clh) \otimes   \clk(\clh))$ such that 
$$ \Sigma^2_m\tau(b)=(\pi \circ \overline{\Sigma^2_m\tau})(b) \,\, \mbox{ for all } b \in \Sigma^2_mA.$$
Let $P_0=1\otimes 1\otimes p$.  Then we have 
$$P_0\sim_{MVN} 1-P_0\sim_{MVN} 1.$$ 
Thus, there exists a isometry $T_0\in M(C(Y) \otimes \clk(\clh) \otimes   \clk(\clh))$ such that 
$$ T_0T_0^*=P_0.$$
This induces the following isomorphism:
$$ T_0\bm{\cdot}T_0^*: M(C(Y) \otimes \clk(\clh) \otimes   \clk(\clh))\rightarrow P_0M(C(Y) \otimes \clk(\clh) \otimes   \clk(\clh))P_0;\quad  a \mapsto T_0aT_0^*.$$

Hence we have 
\begin{center}
	\begin{tikzcd}[
		arrows={-{Stealth[length=5pt,width=5pt]}}, % bigger arrowheads
		line width=1.5pt,                         % thicker shafts
		column sep=huge, row sep=huge             % optional: more spacing
		]
		A \arrow[r, "\Sigma^2_m\tau"] \arrow[dr, swap, "\overline{\Sigma^2_m\tau}"] 
		&Q(C(Y) \otimes \clk(\clh) \otimes   \clk(\clh))  \arrow[r, "i"]  & {[P_0]\,Q(C(Y) \otimes \clk(\clh) \otimes \clk(\clh))[P_0]}\\
		&M(C(Y) \otimes \clk(\clh) \otimes   \clk(\clh))\arrow[r, "T_0\,\bm{\cdot}\,T_0^*"]  
		\arrow[u, "\pi"]  &P_0M(C(Y) \otimes \clk(\clh) \otimes   \clk(\clh))P_0	\arrow[u, "\pi"]
	\end{tikzcd}
\end{center}
Note that $P_0M(C(Y) \otimes \clk(\clh) \otimes   \clk(\clh))P_0=M(C(Y) \otimes \clk(\clh))\otimes p$. Using this and the commutatice diagram, we  get a $*$-homomorphism 
$$\tau^{\prime}: A \rightarrow M(C(Y) \otimes \clk(\clh); \, \tau^{\prime}(a)=\pi\circ \tau.$$ This proves that $[\tau]_{su}=0$. 
\qed \\
We will show that $\Psi_m$ is an isomorphism. To do this, we need to establish some results. Take a homogeneous extension $\lambda: \Sigma_m^2A \rightarrow Q\big(C(Y)\otimes \clk(\clh) \otimes \clk(\clh)\big)$.  Let $\lambda(1\otimes (S^*)^m)=V$. Since $\lambda$ is homoeneous, hence norm full extension, it follows from Proposition $2.5$, and Corollary $2.7$ that there exists a unitary $U\in M\big(C(Y)\otimes \clk(\clh) \otimes \clk(\clh)\big)$ such that 
$[U]V[U^*]=[1\otimes 1 \otimes (S^*)^m]$.  Since $[\lambda]_{su}=[[U]\lambda[U^*]]_{su}$,  one can, without loss of generality assume that $V=[1\otimes 1 \otimes (S^*)^m]$. Hence we have 
 $$\lambda(1\otimes p_{<m})=[1\otimes 1 \otimes p_{<m}].$$
Define 
 $$\nu: M_m(\bbc)\rightarrow  Q\big(C(Y)\otimes \clk(\clh) \otimes \clk(\clh)\big); \, p_{ij}\mapsto \lambda(1\otimes p_{ij}) \mbox{ for all } 1\leq i,j \leq m-1.$$
 Then $\nu \in \Delta_m^2$. Moreover,  if we define 
 $$\tilde{\tau}:A \rightarrow Q\big(C(Y)\otimes \clk(\clh) \otimes \clk(\clh)\big); \, a \mapsto \lambda(a\otimes p_{<m}) \mbox{ for } a \in A, $$ then we get 
 $$\lambda=\Sigma^2_m(\tilde{\tau}, \nu, [1\otimes 1 \otimes (S^*)^m]).$$ 
\bppsn \label{order}
One has 
$$ \mbox{Ord}\,([\nu])\leq \mbox{Ord}\,([\Sigma^2_m(\tilde{\tau}, \nu, [1\otimes 1 \otimes (S^*)^m])]_{su}).$$
In particular, if  $[\nu]$ is of infinite order, then so is $[\Sigma^2_m(\tilde{\tau}, \nu, [1\otimes 1 \otimes (S^*)^m])]_{su}$. 
\eppsn 
\prf  Let $ \mbox{Ord}\,([\Sigma^2_m(\tilde{\tau}, \nu, [1\otimes 1 \otimes (S^*)^m])]_{su})=r$. Fix $r$ isometries $s_1,s_2,\cdots s_r \in Q(C(Y) \otimes \clk(\clh) \otimes   \clk(\clh))$ such that $\sum_{j=1}^r\, s_js_j^*=1$. Define 
$$\phi: \Sigma_m^2A \rightarrow Q(C(Y) \otimes \clk(\clh) \otimes   \clk(\clh)); \quad \phi(b)=\sum_{i=1}^r s_i \Sigma_m^2 (\tau, \nu)(b)s_i^*.$$ 
 Then $[\phi]_{su}=r[\Sigma_m^2 (\tau, \nu)]_{su}$, and hence we have the following commutative diagram:
	\begin{center}
	\begin{tikzcd}[
		arrows={-{Stealth[length=9pt,width=9pt]}}, % bigger arrowheads
		line width=1pt,                         % thicker shafts
		column sep=huge, row sep=huge             % optional: more spacing
		]
	\Sigma_m^2A\arrow[r, "\tilde{\phi}"] \arrow[dr, swap, "\phi"]
		&M(C(Y) \otimes \clk(\clh) \otimes   \clk(\clh)) \arrow[d, "\pi"] \\
		&Q(C(Y) \otimes \clk(\clh) \otimes   \clk(\clh))
	\end{tikzcd}
\end{center}

By restricting it to the $C^*$-subalgebra of $\Sigma_m^2A$ generated by $\{1\otimes p_{ij}:1\leq i,j \leq m\}$, we get the  following commutative diagram:
	\begin{center}
	\begin{tikzcd}[
		arrows={-{Stealth[length=9pt,width=9pt]}}, % bigger arrowheads
		line width=1.5pt,                         % thicker shafts
		column sep=huge, row sep=huge             % optional: more spacing
		]
	M_m(\bbc) \arrow[r, "\tilde{\phi}_{|M_m(\bbc)}"] \arrow[dr, swap, "\nu"]
		&M(C(Y) \otimes \clk(\clh) \otimes   \clk(\clh)) \arrow[d, "\pi"] \\
		&Q(C(Y) \otimes \clk(\clh) \otimes   \clk(\clh))
	\end{tikzcd}
\end{center}

This proves that $r[\nu]=0$, hence the claim. 

\qed 
 \blmma \label{main lemma} 
$[\nu]=0$ if and only if  $[\Sigma^2_m(\tilde{\tau}, \nu, [1\otimes 1 \otimes (S^*)^m])]_{su}\in 	\Psi_m(\operatorname{Ext}_{\mathrm{PPV}}(Y, A))$.
\elmma  
\prf Let $[\Sigma^2_m(\tilde{\tau}, \nu, [1\otimes 1 \otimes (S^*)^m])]_{su}\in 	\Psi_m(\operatorname{Ext}_{\mathrm{PPV}}(Y, A))$. Then there exists a homogeneous extension  $\tau^{\prime}$  such that 
$$\Sigma^2_m(\tilde{\tau}, \nu, [1\otimes 1 \otimes (S^*)^m]) \sim_{su} \Sigma_m^2\tau^{\prime}.$$
Therefore, there exists a unitary $U \in M(C(Y) \otimes \clk(\clh) \otimes   \clk(\clh))$ such that 
$$ [U] \Sigma^2_m(\tilde{\tau}, \nu, [1\otimes 1 \otimes (S^*)^m])(b) [U]^*= \Sigma_m^2\tau^{\prime}(b) \quad \mbox{ for all } b \in  \Sigma_m^2A.$$
Hence we have 
$$ [U] \nu(p_{ij})[U]^*= [U] \Sigma_m^2 (\tau, \nu)(1 \otimes p_{ij}) [U]^*= \Sigma_m^2\tau^{\prime}(1 \otimes p_{ij})=1\otimes 1 \otimes p_{ij}. $$
This shows that $\nu \sim_{su} 0$, hence $[\nu]=0$. To show the forward direction, we assume that  $[\nu]=0$.  
Let $P_{<m}=1\otimes 1 \otimes p_{<m}$. Since 
$$\nu(1)=[P_{<m}], \mbox{ and } 1\otimes P_{<m}\sim_{MVN}1 \sim_{MVN}(1-P_{<m}),$$
there exists a isometry $T_{<m}\in M(C(Y) \otimes \clk(\clh) \otimes   \clk(\clh))$ such that 
$$ T_{<m}T_{<m}^*=P_{<m}.$$
This induces the following isomorphism:
$$ T_{<m}\bm{\cdot}T_{<m}^*: M(C(Y) \otimes \clk(\clh) \otimes   \clk(\clh))\rightarrow P_{<m}M(C(Y) \otimes \clk(\clh) \otimes   \clk(\clh))P_{<m};\quad  a \mapsto  T_{<m}a T_{<m}^*.$$
Hence we have 
	\begin{center}
	\begin{tikzcd}[
		arrows={-{Stealth[length=5pt,width=5pt]}}, % bigger arrowheads
		line width=1.5pt,                         % thicker shafts
column sep=huge, row sep=huge             % optional: more spacing
		]
		M_m(\bbc) \arrow[r, "\nu"] \arrow[dr, swap, "\tilde{\nu}"] 
		&Q(C(Y) \otimes \clk(\clh) \otimes   \clk(\clh))  \arrow[r, "i"]  & {[P_{<m}]Q(C(Y) \otimes \clk(\clh) \otimes   \clk(\clh))[P_{<m}]} \\
		&M(C(Y) \otimes \clk(\clh) \otimes   \clk(\clh))\arrow[r, "T_{<m}\,\bm{\cdot}\,T_{<m}^*"]  
		\arrow[u, "\pi"]  &P_{<m}M(C(Y) \otimes \clk(\clh) \otimes   \clk(\clh))P_{<m}	\arrow[u, "\pi"].
	\end{tikzcd}
\end{center}
Therefore, there exists a unitary $W\in P_{<m}M(C(Y) \otimes \clk(\clh) \otimes   \clk(\clh))P_{<m}$ such that 
$$[W]\nu(p_{ij})[W^*]=[1\otimes 1 \otimes p_{ij}], \mbox{ for all } 0\leq i,j \leq m-1.$$
Define $V=\sum_{j=0}^\infty\, W(1\otimes 1\otimes S^m)^j$. Then $V$ is a unitary in $M(C(Y) \otimes \clk(\clh) \otimes   \clk(\clh))$, and 
$$ [V][1\otimes 1\otimes S^m]=[1\otimes 1\otimes S^m][V], \mbox{ and } [V]\nu(p_{ij})[V^*]=[1\otimes 1 \otimes p_{ij}], \mbox{ for all } 0\leq i,j \leq m-1.$$
Define  the extension $\Gamma=V[\Sigma_m^2 (\tau, \nu)]V^*$. Then $[\Gamma]_{su}=[\Sigma_m^2 (\tau, \nu)]_{su}$. Moreover, the extension $\Gamma$ induces a homomorphism 
$\gamma:A \rightarrow Q(C(Y) \otimes \clk(\clh)$ such that 
$$\Gamma(a\otimes p)=\gamma(a) \otimes p, \mbox{ for all } a \in A.$$
It is not difficult to verify that $\Gamma=\Sigma^2_m(\gamma)$. Therefore, we have 
$$ [\Sigma_m^2 (\tau, \nu)]_{su}=[\Gamma]_{su} \in 	\Psi_m(\operatorname{Ext}_{\mathrm{PPV}}(Y, A)).$$
\qed 

\blmma \label{a1} 
	Let $Y$ be a finite-dimensional compact metric space. Suppose that the groups 
$K_0(C(Y))$ and $K_1(C(Y))$ are finitely generated free abelian groups. Then, for $r\in \bbn_0$, we have $$ \operatorname{Ext}_{\mathrm{PPV}}^{\mathrm{Tor}}(Y, \Sigma_m^{2r}\bbc)=0.$$
\elmma

\prf We prove the claim by induction on $r$. The case $r=0$ is clear, since  
\[
\operatorname{Ext}_{\mathrm{PPV}}(Y,\mathbb{C}) = 0. 
\]
Assume that the claim holds for $r-1$. Let  
$
[\psi]_{su} \in \operatorname{Ext}_{\mathrm{PPV}}^{\mathrm{Tor}}(Y,\Sigma_m^{2r}\mathbb{C}).
$
Then $$[\psi]_{su} = [\Sigma_m^2\bigl(\tau, \, \nu,  [1 \otimes 1 \otimes S^m] \bigr)]_{su},$$
for some $\nu \in \Delta_m^2$ and some homogeneous extension  $\tau$ of   $\Sigma_m^{2r-2}\bbc$  by $C(Y) \otimes \clk$. 
Since
\[
\operatorname{Ext}(M_m(\mathbb{C}), C(Y)) = KK^1(M_m(\mathbb{C}), C(Y))
\]
is a torsion-free abelian group by the K\"unneth Theorem (see Theorem~23.1.2, p.~234 of \cite{Bla-1998aa}), it follows that
\[
\operatorname{Ord}([\nu]) \in \{0, \infty\}.
\]
By Proposition~\ref{order}, we have
\[
\operatorname{Ord}([\nu]) \leq \operatorname{Ord}\!\left([\Sigma_m^2\bigl(\tau, \, \nu,  [1 \otimes 1 \otimes S^m] \bigr)]_{su}\right) < \infty.
\]
Hence $\operatorname{Ord}([\nu]) =0$, and therefore by Lemma~\ref{main}
\[
\psi \in \Psi_m\bigl(\operatorname{Ext}_{\mathrm{PPV}}(Y,\Sigma_m^{2r-2}\mathbb{C})\bigr).
\]
The claim now follows from Lemma \ref{injective hom} and the induction hypothesis. 
\qed

\blmma \label{orderzero} 
	Let $Y$ be a finite-dimensional compact metric space. Suppose that the groups 
$K_0(C(Y))$ and $K_1(C(Y))$ are finitely generated free abelian groups. Let $\nu \in \Delta_m^{\ell}$. Then we have the following.
$$ \mbox{Ord}\,([\nu])=0.$$
\elmma 
\prf By  the K\"unneth Theorem (see Theorem 23.1.2, p.~234, \cite{Bla-1998aa}), 
\[
\operatorname{Ext}(M_m(\mathbb{C}), C(Y)) \cong KK^1(M_m(\mathbb{C}), C(Y))
\]
is a torsion-free abelian group. Hence  we get 
\[
\operatorname{Ord}([\nu]) \in \{0,\infty\}.
\]
Assume that $\operatorname{Ord}([\nu]) = \infty$. Then, taking $A = \mathbb{C}$ and using Proposition~\ref{order} and Lemma~\ref{main lemma}, we obtain
\[
\operatorname{Ord}([\Sigma_m^2 (0, \nu)]) = \infty, 
\quad \text{and} \quad 
\Sigma_m^2 (0, \nu) \notin \Psi_m(\operatorname{Ext}_{\mathrm{PPV}}(Y, \mathbb{C})).
\]
Fix $r \in \bbn_0$. Consider the map
\[
\Psi_m^r: \operatorname{Ext}_{\mathrm{PPV}}(Y, \Sigma_m^{2r}A) 
\longrightarrow 
\operatorname{Ext}_{\mathrm{PPV}}(Y, \Sigma_m^{2(r+1)}A), 
\quad 
[\tau]_{su} \mapsto [\Sigma_m^2 \tau]_{su}.
\]
By a similar argument,
\[
\operatorname{Ord}([\Sigma_m^2 (0_r, \nu)]) = \infty, 
\quad \text{and} \quad 
\Sigma_m^2 (0_r, \nu) \notin 
\Psi_m^r\big(\operatorname{Ext}_{\mathrm{PPV}}(Y, \Sigma_m^{2r}\mathbb{C})\big).
\]
Thus, the elements
\[
\left\{ 
\Psi_m^{r} \circ \cdots \circ \Psi_m^{j+1} \Sigma_m^2 (0_j, \nu) 
\;\middle|\; 0 \le j \le r 
\right\}
\]
are $\mathbb{Z}$-linearly independent elements of infinite order in  
$\operatorname{Ext}_{\mathrm{PPV}}(Y, \Sigma_m^{2r}\mathbb{C})$.
By Proposition \ref{injectivekk}, one can identify  
$\operatorname{Ext}_{\mathrm{PPV}}(Y, \Sigma_m^{2r}\mathbb{C})$ as a subgroup of  
$KK^1(C(Y), \Sigma_m^{2r}\mathbb{C})$. Moreover, by the Universal Coefficient Theorem (see Theorem 23.1.1, p.~233, \cite{Bla-1998aa}), we have
\[
KK^1(\Sigma_m^{2r}\mathbb{C}, C(Y))
\cong 
K_0(C(Y)) \oplus K_1(C(Y)) \oplus (\mathbb{Z}/m\mathbb{Z})^{\oplus r}.
\]
Choose $r$ greater than the rank of $K_{*}(C(Y))$. Then  
$\operatorname{Ext}_{\mathrm{PPV}}(Y, \Sigma_m^{2r}\mathbb{C})$ contains at least $r+1$ $\mathbb{Z}$-linearly independent elements of infinite order, which is possible only if
\[
\operatorname{Ext}_{\mathrm{PPV}}^{\mathrm{Tor}}(Y, \Sigma_m^{2r}\mathbb{C}) \neq 0,.
\]
This leads to  a contradiction to Lemma \ref{a1}. 
\qed 
\bthm \label{main}
Let $m \geq 1$, and  let $A$ be a unital, nuclear, separable 
$C^*$-algebra. Let $Y$ be a finite-dimensional compact metric space such that  $K_0(C(Y))$ and  $K_1(C(Y))$ are finitely generated free abelian groups. Then the map
\[
\Psi_m: \operatorname{Ext}_{\mathrm{PPV}}(Y, A) \to \operatorname{Ext}_{\mathrm{PPV}}(Y, \Sigma_m^2 A), \quad [\tau]_u \mapsto [\Sigma_m^2 \tau]_u
\]
is an isomorphism.
\ethm 
\prf  
It suffices to show that $\Psi_m$ is surjective, by Lemma~\ref{injective hom}. 
Let $\lambda$ be a homogeneous extension of $\Sigma_m^2 A$ by $C(Y) \otimes \mathcal{K}$. 
Then
\[
[\lambda]_{su}
=
\big[\Sigma_m^2(\tilde{\tau}, \nu, [1 \otimes 1 \otimes (S^*)^m])\big]_{su}
\]
for some $\nu \in \Delta_m^2$ and some homogeneous extension $\tau$ of $A$ by $C(Y) \otimes \mathcal{K}$. 
By Lemma~\ref{orderzero}, we have $[\nu] = 0$. 
Therefore, Lemma~\ref{main lemma} implies that
\[
\Sigma_m^2(\tilde{\tau}, \nu, [1 \otimes 1 \otimes (S^*)^m])
\in
\Psi_m\big(\operatorname{Ext}_{\mathrm{PPV}}(Y, A)\big).
\]

\bppsn \label{Shortexactseq}
Let 
\[
\xi:\; 0 \longrightarrow C(Y)\otimes \clk \longrightarrow B \xrightarrow{\,\,\alpha\,\, } A \longrightarrow 0
\]
be an extension of $A$ by $C(Y)\otimes \clk$, and let $\tau$ denote the corresponding Busby invariant. Then 
\[
\Sigma^2_m\xi:\; 0 \longrightarrow C(Y)\otimes \clk  \otimes \clk \longrightarrow \Sigma^2_m B \xrightarrow{\,\, \Sigma_m^2\alpha \,\,} \Sigma^2_m A \longrightarrow 0
\]
is an extension of $\Sigma^2_m A$ by $C(Y)\otimes \clk$ with Busby invariant $\Sigma^2_m \tau$. Moreover, if $\xi$ is homogeneous, then so is $\Sigma^2_m \xi$.
\eppsn
\prf The first part follows by a straightforward verification. For the other part, fix \(y\in Y\) and let
\[
J_y=\ker\big(\mathrm{ev}_y\circ\Sigma^2_m\tau\big).
\]
Since \(\tau\) is homogeneous, \(J_y\) is an ideal of \(\Sigma^2_m A\) with
\(J_y\cap (A\otimes\mathcal K)=\emptyset\). If \(J_y\neq\{0\}\), then there exists \(t\in\mathbb T\) such that
\[
(1\otimes S^m)-t\one \in J_y.
\]
However, we have   
\[
\mathrm{ev}_y\circ\Sigma^2_m\tau(1\otimes S^m)=[S^m] \neq [t\one]
\quad\text{for any }t\in\mathbb T,
\]
so no such nonzero element  lies in $J_y$. Hence $J_y=\{0\}$ for every $y\in Y$, proving the claim.
\qed 
\bppsn \label{extensions} 	Let $Y$ be a finite-dimensional compact metric space. Suppose that the groups 
$K_0(C(Y))$ and $K_1(C(Y))$ are finitely generated free abelian groups. Suppose that 
\[
\mathrm{Ext}_{\mathrm{PPV}}(Y, A)=\{[\xi_i]_{su} : i \in I\},
\]
and let $C_i$ denote the $C^*$-algebra appearing as the middle $C^*$-algebra in the
extension $\xi_i$. Then
\[
\mathrm{Ext}_{\mathrm{PPV}}(Y, \Sigma_m^2 A)
= \{[\Sigma_m^2 \xi_i]_{su} : i \in I\}.
\]
Moreover, the collection
\[
\{\Sigma_m^2 C_i : i \in I\}
\]
is precisely the set of all $C^*$-algebras that appear as middle
$C^*$-algebras of homogeneous extensions of $\Sigma_m^2 A$ by
$C(Y)\otimes \clk$, up to isomorphism.
\eppsn
\prf It is a direct consequence of Theorem \ref{main} and Proposition \ref{Shortexactseq}.
\qed

\bppsn  \label{complete invariant} 
Suppose that $K_0$ is a complete invariant for the middle $C^*$-algebras appearing in the extensions of $\mathrm{Ext}_{\mathrm{PPV}}(Y, A)$. Then $K_0$ is also a complete invariant for the middle $C^*$-algebras appearing in the extensions of $\operatorname{Ext}_{\mathrm{PPV}}(Y, \Sigma^2_m A)$, for all $ m \geq 1$.
\eppsn 
\prf 
It follows from Proposition \ref{Shortexactseq} and Theorem \ref{K-groups for m-QDS}.
\qed  \\
We now establish the main result.
\bthm 
For $1 \leq k \leq 2n$, we have
\[
B_k^{2n+1} \;\cong\; 
\begin{cases}
	\Sigma^{2(k-1)} C(\bbbt), & \text{if } 1 \leq k \leq n, \cr
	\Sigma^{2(k-n-1)} \,\Sigma_2^2 \,\Sigma^{2(n-1)} C(\bbbt), & \text{if } n+1 \leq k \leq 2n. \cr 
\end{cases}
\]
In particular, 
\[
C\bigl(SO_q(2n+1)/SO_q(2n-1)\bigr) \;\cong\; \Sigma^{2(n-1)} \,\Sigma_2^2 \,\Sigma^{2(n-1)} C(\bbbt).
\]
\ethm

\prf Let $1 \leq k \leq n$. In this case, it follows from the explicit description of irreducible representations given in \cite{BhuBisSau-2024aa} that the images of the first $k$ generators of $C\bigl(SO_q(2n+1)/SO_q(2n-1)\bigr)$ under the map $\phi_{\omega_k}$ form the standard generators of $C\bigl(S_q^{2k+1}\bigr)$, while the images of the remaining generators are zero.
The claim then follows, since
\[
C\bigl(S_q^{2k+1}\bigr)\cong \Sigma^{2(k-1)} C(\bbbt).
\]
To prove the claim for $k=n+1$, assume that it holds for $1 \leq k \leq n$.
For $m \neq 0$, we have the extension
\[
\xi_m:\; 0 \longrightarrow C(\bbbt)\otimes \clk \longrightarrow \Sigma_m^2 C(\bbbt)
\xrightarrow{\sigma} C(\bbbt) \longrightarrow 0,
\qquad \sigma(S^m)=\bbt.
\]
To get a  trivial homogeneous extension $\xi_0$, take $q \in \Upsilon$, where  $\Upsilon:=\{q \in \bbc: 0<|q| <1, \theta=\frac{1}{\pi}\arg{(q)} \mbox{ is irrational}\}$. Then by \cite{Sau-2019aa}, we have the following trivial homogeneous extension. 
\[ 
\xi_0: 0 \longrightarrow  C(\bbbt)\otimes \clk  \longrightarrow C(U_q(2)/_{\psi}\bbbt) 
\stackrel{\vartheta}{\longrightarrow} C(\bbbt) \longrightarrow 0.
\]
Denote $C(U_q(2)/_{\psi}\bbbt)$ by $\Sigma^2_0C(\bbbt)$.  Then we have 
$$ K_0(\Sigma^2_0C(\bbbt))=\begin{cases}
\bbz\oplus \bbz/m\bbz& \mbox{  if } m \in \bbz\setminus \{0\}, \cr 
\bbz \oplus \bbz  & \mbox{ if } m=0. \cr
\end{cases}$$
Using Lemma $3.4$ in \cite{Sau-2019aa}, we get 
\[
\mathrm{Ext}_{\mathrm{PPV}}(\bbbt, C(\bbbt))=\{[\xi_m]_{su}: m\in \bbz\}.
\]
Since $ \Sigma_m^2 C(\bbbt)= \Sigma_{-m}^2 C(\bbbt)$, it follows that $K_0$ is a complete invariant for the middle $C^*$-algebras appearing
in the extensions of $\operatorname{Ext}_{\mathrm{PPV}}(\bbbt, C(\bbbt))$. Then 
from Proposition~\ref{complete invariant}, one can conclude that $K_0$ is also a complete invariant for
the middle $C^*$-algebras appearing in the extensions of
$\operatorname{Ext}_{\mathrm{PPV}}(\bbbt, B_n)$.
By Proposition~\ref{extensions}, we obtain
\[
\mathrm{Ext}_{\mathrm{PPV}}(\bbbt, B_n^{2n+1})
= \mathrm{Ext}_{\mathrm{PPV}}(\bbbt, \Sigma^{2(n-1)} C(\bbbt))
= \{[\Sigma^{2(n-1)} \xi_m]: m\in \bbz\}.
\]
Therefore, the middle $C^*$-algebras appearing in the extensions of
$\operatorname{Ext}_{\mathrm{PPV}}(\bbbt, B_n)$ are given by 
$$\{ \Sigma^2_m\Sigma^{2(n-1)}C(\bbbt): m \in \bbz\}.$$ By Theorem \ref{K-groups for m-QDS},  we have 
$$ K_0(\Sigma^2_m\Sigma^{2(n-1)}C(\bbbt))=\begin{cases}
	\bbz\oplus \bbz/m\bbz& \mbox{  if } m \in \bbz\setminus \{0\}, \cr 
	\bbz \oplus \bbz  & \mbox{ if } m=0. \cr
\end{cases}$$
By the induction hypothesis and Lemma~\ref{lemma-homogeneous}, it follows that
$B_{n+1}^{2n+1}$ appears as a middle $C^*$-algebra in the extensions belonging to
$\mathrm{Ext}_{\mathrm{PPV}}(\bbbt, B_n^{2n+1})$.
Now, using Theorem~\ref{K-groups for m-QDS} and \ref{K-groups-B}, we conclude that
\[ 
B_{n+1}^{2n+1}\cong \Sigma_2^2 \Sigma^{2(n-1)} C(\bbbt).
\]
For the case $k=n+2$, the argument is the same as given above. However, in this case one needs to use Theorem~\ref{main} or Proposition \ref{extensions} to obtain the list of all extensions, and hence all the middle $C^*$-algebras appearing in the extensions of $\operatorname{Ext}_{\mathrm{PPV}}(\bbbt, B_{n+1})$. The remaining cases follow along similar lines. 
\qed 

\bcrlre
	For all $q \in (0,1)$, the $C^*$-algebras $C\bigl(SO_q(2n+1)/SO_q(2n-1)\bigr)$ are isomorphic.

\ecrlre

%%%%%%%%%%%%%%%%%%%%%%%%%%%%%%%%%
%%%%%%%%%%%    Another action   %%%%%%%%%%%%
%%%%%%%%%%%%%%%%%%%%%%%%%%%%%%%%%

\bigskip

\bigskip

\noindent{\sc Bipul Saurabh} (\texttt{bipul.saurabh@iitgn.ac.in},  \texttt{saurabhbipul2@gmail.com})\\
{\footnotesize Department of Mathematics,  Indian Institute of Technology, Gandhinagar, Palaj, Gandhinagar 382055, India}

\end{document}